\documentclass[12pt]{article}
\usepackage{amssymb}

\usepackage[]{amsfonts}

\def\underset#1#2{\mathrel{\mathop{\kern0pt #2}\limits_{#1}}}

\def\overset#1#2{\mathrel{\mathop{\kern0pt #2}\limits^{#1}}}

\def\couleur(#1 #2 #3)
	{
\special{ps::[end]}}

\def\sqr#1#2{{\vcenter{\vbox{\hrule height.#2pt
             \hbox{\vrule width.#2pt height#1pt \kern#1pt
             \vrule width.#2pt}
             \hrule height.#2pt}}}}

\def\square{\mathchoice\sqr64\sqr64\sqr43\sqr23}			

\def\st{\mathinner{\mkern1mu\raise1pt\hbox{.}				
		   \mkern1mu\raise4pt\hbox{.}
		   \mkern1mu\raise1pt\hbox{.}
		 }
         }

\def\bx#1{\setbox1=\hbox{\kern3pt{#1}\kern3pt}				
 \dimen1=\ht1 \advance\dimen1 by 3pt \dimen2=\dp1 \advance\dimen2 by 3pt
 \setbox1=\hbox{\vrule height\dimen1 depth\dimen2\box1\vrule}%
 \setbox1=\vbox{\hrule\box1\hrule}%
 \advance\dimen1 by .4pt \ht1=\dimen1
 \advance\dimen2 by .4pt \dp1=\dimen2 \box1\relax}

\def\k#1{\kern#1em}
\def\vci{\vrule  width.02em height1.47ex depth-.0ex}				
\def\11{{\rm\k{.2}\vci\k{-.37}1}}

\newtheorem{Theorem}{Theorem}[section]
\newtheorem{Definition}[Theorem]{Definition}
\newtheorem{Lemma}[Theorem]{Lemma}
\newtheorem{Corollary}[Theorem]{Corollary}
\newtheorem{Proposition}[Theorem]{Proposition}

\begin{document}
\title{Analysis by discs.}
\author{Eric Amar}
\maketitle
\ \par
\ \par
\tableofcontents
\ \par
\section{Introduction.}
\ \ \ \ \ \ The aim of analysis by discs is to see the possibility to 
reduce certain problems in several variables to problems in one variable.\ 
\par
\ \ \ \ \ \ We shall work in the unit ball ${\Bbb B}$ of ${\Bbb C}^{n}$ 
and a disc in ${\Bbb B}$ is a holomorphic mapping $\phi :\ {\Bbb D}{\longrightarrow}{\Bbb B}$. 
The set of all discs in ${\Bbb B}$ will be denoted by ${\cal{D}}({\Bbb B})$.\ 
\par
\ \ \ \ \ \ The analysis by discs can be roughly stated the following 
way:\ \par
suppose you have a several variables object ${\cal{O}}$ and a property 
${\cal{P}}$ true for all discs $\phi $, for the one variable object $"{\cal{O}}\circ \phi "$; 
is the corresponding property ${\cal{P}}$  true for ${\cal{O}}$ ?\ \par
\ \ \ \ \ \ In order to be more specific with examples, let me recall 
the definitions of Hardy and Bergman spaces in the ball.\ \par

\begin{Definition} Let $k\in {\Bbb N}$ then:\ \par
\ \ \ \ \ \ \ \ \ \ $\displaystyle A_{k}^{p}({\Bbb B}):=\{f\in {\cal{O}}({\Bbb B})\ /\ \left\Vert{f}\right\Vert 
_{A_{k}^{p}({\Bbb B})}^{p}:=\int_{{\Bbb B}}^{}{\left\vert{f(z)}\right\vert 
^{p}(1-\left\vert{z}\right\vert ^{2})^{k}\,d\lambda (z)}<\infty \}$.\ 
\par
{\rm{\ \ \ \ \ \ \ \ \ \ \ \ }}$\displaystyle H^{p}({\Bbb B}):=\{f\in {\cal{O}}({\Bbb B})\ /\ \left\Vert{f}\right\Vert 
_{H^{p}({\Bbb B})}^{p}:={\mathrm{su}}{\mathrm{p}}_{r<1}\int_{\partial 
{\Bbb B}}^{}{\left\vert{f(r\zeta )}\right\vert ^{p}\,d\sigma (\zeta )}<\infty 
\}${\rm{.}}\ \par
\end{Definition}
\ \ \ \ \ \ In order to have a natural problem in $H^{p}({\Bbb B})$, the 
first question is:\ \par
if $f\in H^{p}({\Bbb B})$ in what space of analytic functions in the unit 
disc $f\circ \phi $ is ? If $\phi $ is a "flat" disc, e.g. $\forall \zeta \in {\Bbb D},\ \phi (\zeta ):=(\zeta ,0)\in {\Bbb B}\subset 
{\Bbb C}^{2}$, then the "subordination lemma"~\cite{amBerg78} gives us 
that $f\circ \phi $ is in a Bergman space. Our main result generalizes 
this to any disc in ${\Bbb B}$:\ \par

\begin{Theorem} Let $\phi (z):=(\phi _{1}(z),...,\ \phi _{n}(z))$ be an 
analytic disc in the unit ball ${\Bbb B}$ of  ${\Bbb C}^{n}$ such that 
$\phi (0)=0$ and let $1\leq p\leq \infty ,\ f\in H^{p}({\Bbb B})$ then:\ 
\par
\ \ \ \ \ \ \ \ \ \ \ \ \ $\displaystyle \int_{{\Bbb D}}^{}{\left\vert{f\circ \phi (z)}\right\vert 
^{p}(1-\left\vert{z}\right\vert ^{2})^{n-2}\,d\lambda (z)}\leq C\left\Vert{f}\right\Vert 
_{p}^{p}$.\ \par
\end{Theorem}
\ \ \ \ \ \ This theorem allows us to give an example where the analysis 
by discs is possible for $H^{p}$.\ \par
* ($H^{p}$ Corona) Let ${\cal{O}}$ be $g_{1},\ g_{2}$ in $\displaystyle H^{\infty }({\Bbb B})$ 
and ${\cal{P}}$ the property that\ \par
\ \ \ \ \ \ \ \ \ \ \ \ \ $\forall f\in H^{p}({\Bbb B}),\ \exists f_{1},\ f_{2}\in H^{p}({\Bbb B})$ 
with $f=f_{1}g_{1}+f_{2}g_{2}$\ \par
\ then one can easily prove that $\left\vert{g_{1}}\right\vert ^{2}+\left\vert{g_{2}}\right\vert ^{2}\geq 
\delta >0$ and using the Corona theorem in one variable one has:\ \par
(Cp) $\exists C>0,\ \forall \phi \in {\cal{D}}({\Bbb B}),\ \forall h\in A_{n-2}^{p}({\Bbb 
D})\ \exists h_{1},h_{2}\in A_{n-2}^{p}({\Bbb D})\ $\ \par
such that $h_{1}g_{1}\circ \phi +h_{2}g_{2}\circ \phi =h,\ \left\Vert{h_{j}}\right\Vert 
_{A_{n-2}^{p}({\Bbb D})}\leq C,\ j=1,2$.\ \par
Just take $h_{j}:=\chi _{j}h$ where $\chi _{j}\in H^{\infty }({\Bbb D}),\ g_{1}\circ \phi \chi _{1}+g_{2}\circ 
\phi \chi _{2}=1$. I.e. the property ${\cal{P}}$ is still  true uniformly 
for all the discs.\ \par
\ \ \ \ \ \ The analysis by discs is then:\ \par
\textsl{let }$g_{1},g_{2}\in H^{\infty }({\Bbb B})$\textsl{ such that 
(Cp) is true, if }$f\in H^{p}({\Bbb B})$\textsl{ are there }$\displaystyle f_{1},\ f_{2}\in H^{p}({\Bbb B})$\textsl{ 
with}\ \par
\textsl{\ }$f_{1}g_{1}+f_{2}g_{2}=f$\textsl{ ?}\ \par
\ \ \ \ \ \ The answer is yes, because one can prove that the condition 
(Cp) implies that $\displaystyle \left\vert{g_{1}}\right\vert ^{2}+\left\vert{g_{2}}\right\vert 
^{2}\geq \delta >0$ in ${\Bbb B}$ and using the $H^{p}({\Bbb B})$ corona~\cite{amCor91}, 
we have the answer yes.\ \par
\ \par
\ \ \ \ \ \ So the same kind of question arises naturally for the interpolating 
sequences:\ \par
** ($H^{p}$ interpolation) Let ${\cal{O}}$ be $\displaystyle S=\{a_{k}\}_{k\in {\Bbb N}}\subset {\Bbb B}$ 
a sequence of points in the unit ball ${\Bbb B}\subset {\Bbb C}^{n}$. 
The property ${\cal{P}}$ is that $S$ is interpolating for $H^{p}({\Bbb B})$ 
(see definitions in secton~\ref{ballPN97}).\ \par
\ \ \ \ \ \ Our main theorem gives us a necessary condition on discs passing 
through $S$ for $S$ to be interpolating in $H^{p}({\Bbb B})$:\ \par
\ \ \ \ \ \ \textsl{If }$S$ \textsl{is }$H^{p}({\Bbb B})$ \textsl{interpolating 
and if }$(\sigma ,\phi )$ \textsl{is a disc passing through }$S$,\textsl{ 
i.e. }$\phi (\sigma )=S$\textsl{, then }$\sigma $ \textsl{is "almost" 
interpolating for a weighted Bergman class.}\ \par
\ \ \ \ \ \ But in contraste to the $H^{p}$ Corona case,  we prove:\ \par

\begin{Theorem} \label{ballPN86}For any $p\in [1,\infty [$, there is a 
sequence $S$ in ${\Bbb B}_{2}\subset {\Bbb C}^{2}$ such that all disc 
$(\sigma ,\phi )$  passing through $S$  is such that $\sigma $ is interpolating 
for $A^{p}({\Bbb D})$ with a uniformly bounded constant but $S$ is not 
an interpolating sequence for $H^{p}({\Bbb B})$.\ \par
\end{Theorem}
\ \ \ \ \ \ The hypothesis on the sequence $S$ in the previous theorem 
is stronger than the necessary condition and the fact that nevertheless 
$S$ is not $H^{p}({\Bbb B})$ interpolating proves that the analysis by 
discs is \textsl{not pertinent} in the case of $\displaystyle H^{p}({\Bbb B})$ 
with a finite $p$.\ \par
\ \par
\ \ \ \ \ \ These two questions, Corona and interpolating sequences, remain 
open for $H^{\infty }({\Bbb B})$. Following the question asked by J. Globevnik~\cite{GlobevnikDiscs88}, 
with P. Thomas~\cite{AmTho94} we already studied the  interpolating sequences 
for $H^{\infty }({\Bbb B})$ from this point of view, and we obtain partial 
answer.\ \par
\ \ \ \ \ \ This work is organized the following way:\ \par
in the first section we prove that the image by a disc $\phi $ in the 
ball ${\Bbb B}$ of ${\Bbb C}^{n}$ of the Bergman measure $\displaystyle (1-\left\vert{z}\right\vert ^{2})^{n-2}\,d\lambda (z)$ 
of the unit disc is a Carleson measure in ${\Bbb B}$.\ \par
In the second section we characterize the pairs $(a_{j},\ \alpha _{j})_{j=1,...,N}$ 
of points $a_{j}\in {\Bbb B},\ \alpha _{j}\in {\Bbb D}$ such that there 
is an analytic disc $\phi $ from ${\Bbb D}$ to ${\Bbb B}$ sending $\alpha _{j}$ 
to $a_{j}$. This theorem is a generalization of the Pick-Nevanlinna theorem.\ 
\par
In section 3 we recall the definitions of interpolating sequences for 
Hardy and Bergman spaces and, because we shall need it, we recall that 
they are invariant by the automorphisms of the ball.\ \par
Finally in section 4 we build the counter-example, i.e. we prove theorem~\ref{ballPN86}.\ 
\par
\section{A Carleson measure.}
\ \ \ \ \ \ The aim of this section is to prove that if $\phi :{\Bbb D}{\longrightarrow}{\Bbb B}$ 
is an analytic disc in ${\Bbb B}$  such that $\phi (0)=0$, then:\ \par
(*)            $\displaystyle \forall f\in H^{p}({\Bbb B}),\ \int_{{\Bbb D}}^{}{\left\vert{f\circ 
\phi (z)}\right\vert ^{p}(1-\left\vert{z}\right\vert ^{2})^{n-2}\,d\lambda 
(z)}\leq C\left\Vert{f}\right\Vert _{p}^{p}$.\ \par
\ \ \ \ \ \ Let $\mu $ be the measure in the ball ${\Bbb B}$ image by 
the map $\phi $ of the weighted Lebesgue measure $\lambda $ on the unit 
disc ${\Bbb D}$; this means precisely:\ \par
\ \ \ \ \ \ \ \ \ \ \ \ $\forall f\in {\cal{C}}(\overline{{\Bbb B}}),\ \int_{{\Bbb B}}^{}{f\,d\mu 
}=\int_{{\Bbb D}}^{}{f\circ \phi (1-\left\vert{z}\right\vert ^{2})^{n-2}\,d\lambda 
(z)}$.\ \par
\ \ \ \ \ \ On the other hand if the measure $\mu $ is Carleson in ${\Bbb B}$, 
then:\ \par
\ \ \ \ \ \ \ \ \ \ \ \ $\displaystyle \forall f\in H^{p}({\Bbb B}),\ \int_{{\Bbb B}}^{}{\left\vert{f}\right\vert 
^{p}\,d\mu }\leq \left\Vert{\mu }\right\Vert _{C}\left\Vert{f}\right\Vert 
_{p}^{p}$,\ \par
where $\left\Vert{\mu }\right\Vert _{C}$ is the Carleson norm of $\mu $.\ 
\par
\ \ \ \ \ \ Hence in order to get (*) it suffices to prove that $\mu $ 
is a Carleson measure in ${\Bbb B}$.\ \par
\ \ \ \ \ \ By the Carleson-H\"ormander characterization of these measures, 
it is enough to show:\ \par
\ \ \ \ \ \ \ \ \ \ \ \ $\displaystyle \int_{Q(\zeta ,\ \delta )}^{}{\,d\mu }\leq C\delta ^{2}$.\ 
\par
where $\displaystyle Q(\zeta ,\ \delta ):=\{z\in {\Bbb B},\ \left\vert{1-\zeta 
\cdot \overline{z}}\right\vert <\delta \}$ is the pseudo-ball centered 
in $\zeta \in \partial {\Bbb B}$ and of radius $\delta $.\ \par
\ \ \ \ \ \ By invariance by complex rotation of the map $\phi $, we can 
send $\zeta $ to ${\11}=(1,0,...,0)$  and in fact it remains to show that:\ 
\par
\ \ \ \ \ \ \ \ \ \ \ \ $\displaystyle \int_{Q(1,\ \delta )}^{}{\,d\mu }\leq C\delta ^{2}$, 
with $\displaystyle Q({\11},\ \delta ):=\{z\in {\Bbb B},\ \left\vert{1-z_{1}}\right\vert 
<\delta \}$.\ \par
\ \ \ \ \ \ Finally the definition of $\mu $ leads to show:\ \par
\ \ \ \ \ \ \ \ \ \ \ \ $\displaystyle \int_{Q(1,\ \delta )}^{}{\,d\mu }=\int_{{\Bbb D}}^{}{\chi 
_{d(1,\delta )}(t)(1-\left\vert{t}\right\vert ^{2})^{n-2}\,d\lambda (t)}\leq 
C\delta ^{2}$, \ \par
with $\displaystyle \chi _{d(1,\delta )}(t)$ the indicatrix of the disc 
$d(1,\delta ):=\{z\in {\Bbb D}\ /\ \left\vert{1-z}\right\vert <\delta 
\}$.\ \par
\ \ \ \ \ \ This will be achieved via a series of lemmas.\ \par
\ \par

\begin{Lemma} Let $\phi :\ {\Bbb D}{\longrightarrow}{\Bbb D}$ such that 
$\phi (0)=0$ and suppose $\phi $ inner, i.e. the boundary values $\phi ^{*}$ 
of $\phi $ are of modulus one, then the image by $\phi ^{*}$ of the Lebesgue's 
measure $\,d\sigma $ on ${\Bbb T}$ is the Lebesgue's measure $\,d\sigma $ 
on ${\Bbb T}$.\ \par
\end{Lemma}
\ \ \ \ \ \ Proof:\ \par
let $f\in A({\Bbb D})$ and set $\mu :=\phi ^{*}\sigma $ the image of the 
measure $\sigma $ by $\phi ^{*}$, then:\ \par
\ \ \ \ \ \ \ \ \ \ \ \ $\displaystyle \int_{{\Bbb T}}^{}{f\,d\mu }=\int_{{\Bbb T}}^{}{f\circ 
\phi ^{*}\,d\sigma }=f\circ \phi (0)$,\ \par
because $\sigma $ represents $0$ for $A({\Bbb D})$. But $\phi (0)=0$ hence 
we have:\ \par
\ \ \ \ \ \ \ \ \ \ \ \ $\displaystyle \int_{{\Bbb T}}^{}{f\,d\mu }=f(0)$,\ 
\par
hence $\mu $ also represents $0$. By unicity of the representing measure 
of $0$ on ${\Bbb T}$, we must have $\mu =\sigma $.$\hfill \square$\ \par

\begin{Lemma} \label{ballPN1}Let $\mu $ be a probability measure on $\overline{{\Bbb D}}$ 
which represents $0$ i.e.\ \par
\ \ \ \ \ \ \ $\displaystyle \forall f\in A({\Bbb D}),\ \int_{}^{}{f(z)\,d\mu (z)}=f(0)$;\ 
\par
let $\nu $ be its balay\'ee by the Poisson kernel\ \par
\ \ \ \ \ \ \ \ \ \ \ \ \ $\displaystyle \forall g\in {\cal{C}}({\Bbb T}),\ \int_{{\Bbb T}}^{}{g(\zeta 
)\,d\nu (\zeta )}:=\int_{}^{}{\tilde g(z)\,d\mu (z)}$,\ \par
where $\displaystyle \tilde g(z):=\int_{{\Bbb T}}^{}{P(\zeta ;z)g(\zeta )\,d\sigma 
(\zeta )}$ is the harmonic extension of $g$ in ${\Bbb D}$.\ \par
\ \ \ \ \ \ Then $\,d\nu =\,d\sigma $.\ \par
\end{Lemma}
\ \ \ \ \ \ Proof:\ \par
let $f\in A({\Bbb D})$ then $\tilde f_{\mid {\Bbb T}}=f$ because $f$ is 
harmonic in ${\Bbb D}$, hence by definition of $\nu $ we  get:\ \par
\ \ \ \ \ \ \ \ \ \ \ \ $\displaystyle \int_{{\Bbb T}}^{}{f(\zeta )\,d\nu (\zeta )}=\int_{}^{}{f(z)\,d\mu 
(z)}=f(0)$;\ \par
hence $\nu $ is a representing measure of $0$ supported by ${\Bbb T}$ 
and by unicity of the representing measure of $0$ on ${\Bbb T}$ we have 
$\nu $ is the normalized Lebesgue measure $\sigma $ on ${\Bbb T}$.$\hfill 
\square$\ \par

\begin{Lemma} \label{ballPN4}Let $I:={\Bbb T}\cap d(1,\delta )$ and $\tilde I(z)$ 
the harmonic extension of the indicatrix  $\chi _{I}$ of $I$, then the 
set\ \par
\ \ \ \ \ \ \ \ \ \ \ \ \ $A_{\delta }:=\{z\in {\Bbb D}\ /\ \tilde I(z)>\frac{\pi }{4}\}$\ 
\par
contains the disc $d(1,\delta )$.\ \par
\end{Lemma}
\ \ \ \ \ \ Proof:\ \par
we shall do it in the half plane, the computations being easier. The Poisson 
kernel is $\displaystyle P(z,t)=\frac{2y}{(x-t)^{2}+y^{2}}$ with $z=x+iy$; 
let $I=]-\delta ,\ \delta [$  then\ \par
\ \ \ \ \ \ \ \ \ \ \ \ $\displaystyle \tilde I(z)=\int_{-\delta }^{\delta }{\frac{2y}{(x-t)^{2}+y^{2}}\,dt}=\int_{\frac{x-\delta 
}{y}}^{\frac{x+\delta }{y}}{\frac{\,du}{1+u^{2}}}$,\ \par
by the change $\displaystyle u=\frac{t-x}{y}$, so we get:\ \par
\ \ \ \ \ \ \ \ \ \ \ \ $\displaystyle \tilde I(z)={\mathrm{Arctan}}\left({\frac{x+\delta }{y}}\right) 
-{\mathrm{Arctan}}\left({\frac{x-\delta }{y}}\right) $,\ \par
hence if $0\leq x<\delta ,\ 0<y<\delta $ we get $\tilde I(z)\geq {\mathrm{Arctan}}(1)=\frac{\pi }{4}$.$\hfill 
\square$\ \par

\begin{Lemma} Let $\phi \in H^{\infty }({\Bbb D}),\ \phi (0)=0$; for $0\leq \rho <1$ 
let $\forall z\in {\Bbb D},\ \phi _{\rho }(z):=\phi (\rho z)$, then:\ 
\par
\ \ \ \ \ \ \ \ \ \ \ \ \ $\displaystyle \left\vert{\{\zeta \in {\Bbb T}\ /\ \phi _{\rho }(\zeta 
)\in A_{\delta }\}}\right\vert \leq 2\left\vert{I}\right\vert =4\delta 
$.\ \par
\end{Lemma}
\ \ \ \ \ \ Proof:\ \par
let $\mu _{\rho }$ be the image by $\phi _{\rho }$ of the Lebesgue's measure 
$\sigma $ then we have:\ \par
\ \ \ \ \ \ \ \ \ \ \ \ $\displaystyle \int_{}^{}{f(z)\,d\mu _{\rho }(z)}=\int_{{\Bbb T}}^{}{f\circ 
\phi _{\rho }(\zeta )\,d\sigma (\zeta )}=f\circ \phi _{\rho }(0)=f(0)$.\ 
\par
Hence $\mu $ represents $0$; let $\nu $ be the balay\'ee by the Poisson 
kernel of $\mu $, then by   lemma~\ref{ballPN1} $\nu =\sigma $. Now we 
have:\ \par
\ \ \ \ \ \ \ \ \ \ \ \ $\displaystyle \int_{}^{}{\tilde I(z)\,d\mu _{\rho }(z)}=\int_{{\Bbb T}}^{}{\chi 
_{I}(\zeta )\,d\nu (\zeta )}=\int_{{\Bbb T}}^{}{\chi _{I}(\zeta )\,d\sigma 
(\zeta )}=\left\vert{I}\right\vert =2\delta $.\ \par
On the other hand by definition of $\mu $ we have:\ \par
\ \ \ \ \ \ \ \ \ \ \ \ $\displaystyle \left\vert{I}\right\vert =\int_{}^{}{\tilde I(z)\,d\mu 
_{\rho }(z)}=\int_{{\Bbb T}}^{}{\tilde I\circ \phi _{\rho }(\zeta )\,d\sigma 
(\zeta )}$,\ \par
hence $\displaystyle \int_{\{\tilde I(z)>\pi /4\}}^{}{\frac{\pi }{4}\,d\mu _{\rho 
}(z)}\leq \int_{\{\tilde I(z)>\pi /4\}}^{}{\tilde I(z)\,d\mu _{\rho }(z)}\leq 
\left\vert{I}\right\vert $,\ \par
and $\displaystyle \mu _{\rho }(A_{\delta })\leq \frac{4}{\pi }\left\vert{I}\right\vert 
$.\ \par
But $\displaystyle \mu _{\rho }(A_{\delta })=\int_{}^{}{\chi _{A_{\delta }}(z)\,d\mu 
_{\rho }(z)}=\int_{{\Bbb T}}^{}{\chi _{A_{\delta }}\circ \phi _{\delta 
}(\zeta )\,d\sigma (\zeta )}=\left\vert{\{\zeta \in {\Bbb T}\ /\ \phi 
_{\rho }(\zeta )\in A_{\delta }\}}\right\vert $,\ \par
hence $\displaystyle \left\vert{\{\zeta \in {\Bbb T}\ /\ \phi _{\rho }(\zeta 
)\in A_{\delta }\}}\right\vert \leq \frac{4}{\pi }\left\vert{I}\right\vert 
$.\ \par

\begin{Theorem} Let $\phi (z):=(\phi _{1}(z),...,\ \phi _{n}(z))$ be an 
analytic disc in the unit ball ${\Bbb B}$ of  ${\Bbb C}^{n}$  such that 
$\phi (0)=0$, then the image $\mu $ in ${\Bbb B}$ of the Bergman measure 
$(1-\left\vert{z}\right\vert ^{2})^{n-2}\,d\lambda (z)$ on ${\Bbb D}$ 
is a Carleson measure in ${\Bbb B}$.\ \par
\end{Theorem}
\ \ \ \ \ \ Proof:\ \par
it is enough to prove that, $\mu (Q({\11},\delta ))\lesssim \delta ^{2}$ 
by rotation, with $\displaystyle {\11}:=(1,0,...,0)\in \partial {\Bbb B}$ 
and\ \par
\ \ \ \ \ \ \ $Q({\11},\delta ):=\{z=(z_{1},...,z_{n})\in {\Bbb B}\ /\ \left\vert{1-z_{1}}\right\vert 
<\delta \}$.\ \par
So $\displaystyle \int_{}^{}{\chi _{Q}(w)\,d\mu (w)}=\int_{{\Bbb D}}^{}{\chi 
_{Q}(\phi _{1}(z),...,\phi _{n}(z))\,d\lambda (z)}\leq \int_{{\Bbb D}}^{}{\chi 
_{A_{\delta }}(\phi _{1}(z))\,d\lambda (z)}$,\ \par
because $Q$ involves only $\phi _{1}$ and $A_{\delta }$ contains the disc 
$d(1,\delta )$ by lemma~\ref{ballPN4}.\ \par
\ \ \ \ \ \ But $\phi _{1}(0)=0$ which means by Schwarz lemma that $\phi _{1}(z)=zg(z)$ 
still with $g:{\Bbb D}{\longrightarrow}{\Bbb D}$, hence in order to have 
that $\phi _{1}(z)\in A_{\delta }$ we need to have:\ \par
\ \ \ \ \ \ \ \ \ \ \ \ \ 
\begin{equation} \left\vert{z}\right\vert \geq 1-\delta \ \label{ballPN2}\end{equation}
\ \par
\ \ \ \ \ \ Now applying the previous lemma to $\phi _{1}$ we get:\ \par
\ \ \ \ \ \ \ \ \ \ \ \ $\displaystyle \left\vert{\{\zeta \in {\Bbb T}\ /\ \phi _{\rho }(\zeta 
)\in A_{\delta }\}}\right\vert \leq \frac{8}{\pi }\delta $.\ \par
We can compute the area of the set $\{z\in {\Bbb D}\ /\ \phi _{1}(z)\in A_{\delta }\}$ 
by cutting it with circles of radius $\rho $ with $1-\delta \leq \rho <1$ 
and we get:\ \par
\ \ \ \ \ \ \ \ \ \ \ \ \ \ \ \ \ \ $\lambda (\{z\in {\Bbb D}\ /\ \phi _{1}(z)\in A_{\delta }\})\leq \left\vert{\{\zeta 
\in {\Bbb T}\ /\ \phi _{\rho }(\zeta )\in A_{\delta }\}}\right\vert \times 
\left\vert{\{1-\delta \leq \rho <1\}}\right\vert \leq \frac{8}{\pi }\delta 
^{2}$,\ \par
and the theorem.$\hfill \square$\ \par
\ \ \ \ \ \ As pointed out by J. Bruna, the one dimensional part of the 
previous proof can  be sligthly shortened by use of the Littlewood's subordination 
lemma(~\cite{Duren70}, p. 10).\ \par

\begin{Corollary} Let $\phi (z):=(\phi _{1}(z),...,\ \phi _{n}(z))$ be 
an analytic disc in the unit ball ${\Bbb B}$ of  ${\Bbb C}^{n}$  such 
that $\phi (0)=0$ and let $1\leq p\leq \infty ,\ f\in H^{p}({\Bbb B})$ 
then:\ \par
\ \ \ \ \ \ \ \ \ \ \ \ \ $\displaystyle \int_{{\Bbb D}}^{}{\left\vert{f\circ \phi (z)}\right\vert 
^{p}(1-\left\vert{z}\right\vert ^{2})^{n-2}\,d\lambda (z)}\leq C\left\Vert{f}\right\Vert 
_{p}^{p}$.\ \par
\end{Corollary}
\ \ \ \ \ \ Proof:\ \par
the image is a Carleson measure and we can apply the Carleson-H\"ormander 
 inequality.$\hfill \square$\ \par
\ \ \ \ \ \ At this point one can ask for the converse of the previous 
corollary:\ \par
suppose that $f\in H^{\infty }({\Bbb B})$ and ${\mathrm{su}}{\mathrm{p}}_{\phi \in {\cal{D}}({\Bbb B})}\left\Vert{f\circ 
\phi }\right\Vert _{A_{n-2}^{p}({\Bbb D})}\leq C$ is it true that $\left\Vert{f}\right\Vert _{H^{p}({\Bbb B})}\leq KC$?\ 
\par
\ \ \ \ \ \ The answer is no and the following counter-example was suggested 
by A. Borichev.\ \par
\ \ \ \ \ \ In the unit ball ${\Bbb B}$ of ${\Bbb C}^{2}$ take a inner 
function $f$ such that $f(0)=0$. This means that $f\in H^{\infty }({\Bbb B})$ 
and $\left\vert{f^{*}}\right\vert =1,\ a.e.\ on\ \partial {\Bbb B}$, where 
$f^{*}$ are the boundary values of $f$. We have that $\forall k\in {\Bbb N},\ \left\Vert{f^{k}}\right\Vert _{p}=1$; 
on the other hand, $f\circ \phi $ is still bounded by $1$ in ${\Bbb D}$ 
and $f\circ \phi (0)=0$  hence by Schwarz lemma we have $\left\vert{f\circ \phi (z)}\right\vert \leq \left\vert{z}\right\vert 
$.\ \par
\ \ \ \ \ \ This implies that $\displaystyle \left\Vert{f^{k}\circ \phi }\right\Vert _{A^{p}({\Bbb D})}^{p}\leq 
\int_{{\Bbb D}}^{}{\left\vert{z}\right\vert ^{kp}\,d\lambda (z)}=\frac{1}{kp+2}$; 
hence\ \par
\ \ \ \ \ \ \ \ \ \ \ \ $\displaystyle {\mathrm{su}}{\mathrm{p}}_{\phi \in {\cal{D}}({\Bbb B})}\left\Vert{f^{k}\circ 
\phi }\right\Vert _{A^{p}({\Bbb D})}\leq \left({\frac{1}{kp+2}}\right) 
^{1/p}{\longrightarrow}0$ when $k{\longrightarrow}\infty $,\ \par
hence the constant $K$ cannot be finite.\ \par
\section{A Pick-Nevanlinna theorem.}
\ \ \ \ \ Modifying the proof of the classical theorem of Pick and Nevanlinna 
by Nagy and Foias, we get:\ \par

\begin{Theorem} Let $\sigma =\{\alpha _{k}\}_{k=1,...,N}\subset {\Bbb D}$ 
and $v=\{v_{k}\}_{k=1,...,N}\subset {\Bbb B}$. There is a function $f:{\Bbb D}{\longrightarrow}{\Bbb B}$ 
such that $\forall k=1,...,N,\ f(\alpha _{k})=v_{k}$  iff:\ \par
\ \ \ \ \ \ \ \ \ \ \ $\displaystyle \forall k=1,...,N,\ \forall h_{k}\in {\Bbb C},\ \sum_{k,l=1}^{N}{\overline{h}_{k}h_{l}\frac{1-\overline{v}_{k}\cdot 
v_{l}}{1-\overline{z}_{k}z_{l}}}\geq 0$.\ \par
\end{Theorem}
\ \ \ \ \ \ If ${\Bbb B}={\Bbb D}$ the unit disc, this is the Pick-Nevanlinna 
theorem.\ \par
\ \par
\ \ \ \ \ \ Let $\alpha \in {\Bbb D}$ and $\displaystyle k_{\alpha }(z)=\frac{1}{1-\overline{\alpha }z}$ 
the Cauchy kernel associated to it and:\ \par
\ \ \ \ \ \ \ \ \ \ \ \ $E_{\sigma }:={\mathrm{span}}\{k_{\alpha },\ \alpha \in \sigma \}\subset 
H^{2}({\Bbb D})$.\ \par
\ \ \ \ \ \ To a function $f\ \in H^{\infty }({\Bbb D})$ let us associate 
the operator:\ \par
\ \ \ \ \ \ \ \ \ \ \ \ $\displaystyle \forall h\in E_{\sigma },\ \pi _{\sigma }(f)h:=P_{E_{\sigma 
}}(\overline{f}h)\in H^{2}({\Bbb D})$.\ \par
This is clearly a anti-holomorphic representation of $H^{\infty }({\Bbb D})$ 
in ${\cal{L}}(H^{2}({\Bbb D}))$ and we have\ \par
\ \ \ \ \ \ \ \ \ \ \ \ $\forall \alpha \in \sigma ,\ \pi _{\sigma }(f)k_{\alpha }=\overline{f}(\alpha 
)k_{\alpha }$.\ \par
Now to $f=(f_{1},...,f_{n})\ :\ {\Bbb D}{\longrightarrow}{\Bbb B}$ holomorphic, 
let us associate the operator:\ \par
\ \ \ \ \ \ \ \ \ \ \ \ \ $\forall h\in E_{\sigma },\ \pi _{\sigma }(f)h:=(\pi _{\sigma }(f_{1})h,...,\pi 
_{\sigma }(f_{n})h)\in {\cal{H}}:=H^{2}({\Bbb D})\oplus \cdot \cdot \cdot 
\oplus H^{2}({\Bbb D})$.\ \par

\begin{Lemma} If $f\ :\ {\Bbb D}{\longrightarrow}{\Bbb B}$ then $\left\Vert{\pi _{\sigma }(f)}\right\Vert _{op}\leq 1$.\ 
\par
\end{Lemma}
\ \ \ \ \ \ Proof:\ \par
we have\ \par
\ \ \ \ \ \ $\displaystyle \pi _{\sigma }(f)h=P_{E_{\sigma }\oplus \cdot \cdot \cdot 
\oplus E_{\sigma }}(\overline{f}_{1}h,...,\overline{f}_{n}h)$;\ \par
now the projection is a contraction and we have\ \par
\ \ \ \ \ \ \ $\displaystyle \left\Vert{(\overline{f}_{1}h,...,\overline{f}_{n}h)}\right\Vert 
_{L^{2}({\Bbb T})\oplus \cdot \cdot \cdot \oplus L^{2}({\Bbb T})}^{2}=\sum_{j=1}^{n}{\int_{{\Bbb 
T}}^{}{\left\vert{\overline{f}_{j}(\zeta )h(\zeta )}\right\vert ^{2}\,d\sigma 
(\zeta )}}$,\ \par
hence\ \par
\ \ \ \ \ \ $\displaystyle \left\Vert{(\overline{f}_{1}h,...,\overline{f}_{n}h)}\right\Vert 
_{L^{2}({\Bbb T})\oplus \cdot \cdot \cdot \oplus L^{2}({\Bbb T})}^{2}=\int_{{\Bbb 
T}}^{}{\sum_{j=1}^{n}{\left\vert{\overline{f}_{j}(\zeta )}\right\vert 
^{2}}\left\vert{h(\zeta )}\right\vert ^{2}\,d\sigma (\zeta )}\leq \left\Vert{h}\right\Vert 
_{L^{2}({\Bbb T})}^{2}$.\ \par
So $\displaystyle \left\Vert{\pi _{\sigma }(f)h}\right\Vert ^{2}\leq \left\Vert{(\overline{f}_{1}h,...,\overline{f}_{n}h)}\right\Vert 
_{L^{2}({\Bbb T})\oplus \cdot \cdot \cdot \oplus L^{2}({\Bbb T})}^{2}\leq 
\left\Vert{h}\right\Vert ^{2}$,\ \par
and the lemma. $\hfill \square$\ \par
\ \ \ \ \ \ Proof of the theorem :\ \par
denote by $Z$ the shift operator (multiplication by $z$) on $H^{2}({\Bbb D})$ 
and also the shift on\ \par
\ \ \ \ \ \ \ \ \ ${\cal{H}}:=H^{2}({\Bbb D})\oplus \cdot \cdot \cdot \oplus H^{2}({\Bbb 
D})$,  i.e. $\displaystyle \forall h=(h_{1},...,h_{n})\in {\cal{H}}{\longrightarrow}(Zh)(z)=(zh_{1}(z),...,zh_{n}(z))$.\ 
\par
These operators are isometries and we get easily:\ \par
\ \ \ \ \ \ $\displaystyle \forall f\in H^{\infty }({\Bbb D})^{n},\ \pi _{\sigma }(f)Z^{*}=Z^{*}\pi 
_{\sigma }(f)$,\ \par
because $E_{\sigma }$ is stable by $Z^{*}$ and both diagonalize through 
the $\{k_{\alpha },\ \alpha \in \sigma \}$; in fact $Z^{*}_{\mid E_{\sigma }}=\pi _{\sigma }(z)$.\ 
\par
\ \ \ \ \ \ Now let any $f\in H^{\infty }({\Bbb D})^{n}$ be such that 
$\forall k=1,...,N,\ f(\alpha _{k})=v_{k}$ and set $T:=\pi _{\sigma }(f)\in {\cal{L}}(E_{\sigma },{\cal{H}})$, 
 and suppose that $\left\Vert{T}\right\Vert _{op}\leq 1$ (which is equivalent 
to the condition in the theorem), then we have $Z^{*}T=TZ^{*}$ on $E_{\sigma }$. 
We can apply the lifting of the commutant theorem: there is an operator 
$\tilde T=(\tilde T_{1},...,\tilde T_{n})\in {\cal{L}}(H^{2}({\Bbb D}),\ 
{\cal{H}})$ such that:\ \par
\textsl{(i)}      $\tilde T_{\mid E_{\sigma }}=T$;\ \par
\textsl{(ii)}      $\left\Vert{\tilde T}\right\Vert =\left\Vert{T}\right\Vert $;\ 
\par
\textsl{(iii)}      $Z^{*}\tilde T=\tilde TZ^{*}$.\ \par
As in~\cite{AmMenOp02} \textsl{(iii)} implies that $Z^{*}\tilde T_{j}=\tilde T_{j}Z^{*}$, 
we know that this is the condition for  $\tilde T_{j}$ to be associated 
to  a function $\tilde f_{j}\in H^{\infty }({\Bbb D})$, again as in~\cite{AmMenOp02} 
and $\tilde T$ is the adjoint of the operator $T_{\tilde f}$  of multiplication 
by $\tilde f=(\tilde f_{1},...,\tilde f_{n})$.\ \par
By \textsl{(ii)} we get $\left\Vert{T_{\tilde f}}\right\Vert =\left\Vert{\tilde T}\right\Vert 
\leq 1$ hence $\tilde f\ :\ {\Bbb D}{\longrightarrow}{\Bbb B}$.\ \par
Finally by \textsl{(i)} we get that $\forall k=1,...,N,\ \tilde f(\alpha _{k})=v_{k}$, 
and the theorem.$\hfill \square$\ \par
\ \par
\ \ \ \ \ \ The analogue for the polydisc is straighforward with the usual 
PN theorem:\ \par

\begin{Theorem} Let $\sigma =\{\alpha _{k}\}_{k=1,...,N}\subset {\Bbb D}$ 
and $v=\{v_{k}\}_{k=1,...,N}\subset {\Bbb D}^{n}$. There is a function 
$f:{\Bbb D}{\longrightarrow}{\Bbb D}^{n}$ such that $\forall k=1,...,N,\ f(\alpha _{k})=v_{k}$ 
 iff:\ \par
\ \ \ \ \ \ \ \ \ \ \ $\displaystyle \forall m=1,...,n,\ \forall k=1,...,N,\ \forall h_{k}\in 
{\Bbb C},\ \sum_{k,l=1}^{N}{\overline{h}_{k}h_{l}\frac{1-\overline{v}_{k}^{m}\cdot 
v_{l}^{m}}{1-z_{k}\overline{z}_{l}}}\geq 0$.\ \par
\end{Theorem}
\ \ \ \ \ \ Proof :\ \par
just apply the PN theorem.$\hfill \square$\ \par
\section{Interpolating sequences and automorphisms.\label{ballPN97}}

\begin{Definition} If $\displaystyle S=\{a_{k},\ k\in {\Bbb N}\}$ is a 
sequence of points in ${\Bbb B}$, then we define:\ \par
\ \ \ \ \ \ $\displaystyle \ell _{H}^{p}(S):=\{\lambda =\{\lambda _{k}\}_{k\in {\Bbb 
N}}\ /\ \sum_{k\in {\Bbb N}}^{}{\left\vert{\lambda _{k}}\right\vert ^{p}(1-\left\vert{a_{k}}\right\vert 
^{2})^{n}}=:\left\Vert{\lambda }\right\Vert _{H,p}^{p}<\infty \}$.\ \par
\end{Definition}
This sequence space depends on the sequence $S$. The same for the weighted 
Bergman  classes:\ \par

\begin{Definition} If $\displaystyle S=\{a_{k},\ k\in {\Bbb N}\}$ is a 
sequence of points in ${\Bbb B}$, then we define:\ \par
{\rm{\ \ \ \ \ \ \ \ \ \ \ \ }}$\displaystyle \ell _{A_{l}}^{p}(S):=\{\lambda =\{\lambda _{k}\}_{k\in 
{\Bbb N}}\ /\ \sum_{k\in {\Bbb N}}^{}{\left\vert{\lambda _{k}}\right\vert 
^{p}(1-\left\vert{a_{k}}\right\vert ^{2})^{n+l+1}}=:\left\Vert{\lambda 
}\right\Vert _{A_{l},p}^{p}<\infty \}$.\ \par
\end{Definition}
\ \ \ \ \ \ Now we can define the interpolating sequences for Hardy spaces:\ 
\par

\begin{Definition} We say that $S$ is interpolating for $\displaystyle H^{p}({\Bbb B})$ 
with constant $C$ if:\ \par
\ \ \ \ \ \ \ $\forall \lambda \in \ell _{H}^{p}(S),\ \exists f\in H^{p}({\Bbb B})\ 
/\ \forall k\in {\Bbb N}\ f(a_{k})=\lambda _{k}$ and $\left\Vert{f}\right\Vert _{p}\leq C\left\Vert{\lambda }\right\Vert _{H,p}$.\ 
\par
The set of interpolating sequences for $H^{p}({\Bbb B})$ will be denoted 
by $IH^{p}({\Bbb B})$.\ \par
\end{Definition}
And for the Bergman classes:\ \par

\begin{Definition} We say that $S$ is interpolating for $\displaystyle A_{l}^{p}({\Bbb B})$ 
with constant $C$ if:\ \par
\ \ \ \ \ \ \ $\forall \lambda \in \ell _{A_{l}}^{p}(S),\ \exists f\in A_{l}^{p}({\Bbb 
B})\ /\ \forall k\in {\Bbb N}\ f(a_{k})=\lambda _{k}$ and $\left\Vert{f}\right\Vert _{A_{l}^{p}({\Bbb B})}\leq C\left\Vert{\lambda 
}\right\Vert _{A_{l},p}$.\ \par
The set of interpolating sequences for $A_{l}^{p}({\Bbb B})$ will be denoted 
by $IA_{l}^{p}({\Bbb B})$.\ \par
\end{Definition}
\ \ \ For the unweighted Bergman classes, i.e. $l=0$, we shall set $\displaystyle \ell _{A_{0}}^{p}(S)=\ell _{A}^{p}(S)$\ 
\par
\ \ \ \ \ \ In a paper by Jevtic-Massaneda-Thomas~\cite{JevMasTho96}, 
the authors prove  the following theorem.\ \par

\begin{Theorem} Let $\displaystyle S=\{a_{k}\}_{k\in {\Bbb N}}\subset {\Bbb B}$ 
be an interpolating sequence in $A_{k}^{p}({\Bbb B})$ and $\Phi $ be an 
automorphism of ${\Bbb B}$ then $\Phi (S):=\{\Phi (a_{k})\}_{k\in {\Bbb N}}\subset {\Bbb B}$ 
is still an interpolating sequence in $A_{k}^{p}({\Bbb B})$ with the same 
constant.\ \par
\end{Theorem}
\ \ \ \ \ The following theorem is in W. Rudin's book (~\cite{RudinBall81}, 
section 7.5):\ \par

\begin{Theorem} Let $\psi \in {\mathrm{Aut}}({\Bbb B}),\ a:=\psi ^{-1}(0)$ 
then the linear operator:\ \par
\ \ \ \ \ \ \ \ \ \ \ \ $\displaystyle \forall f\in H^{p}({\Bbb B}),\ Tf(z):=\frac{\left({1-\left\vert{a}\right\vert 
^{2}}\right) ^{n/p}}{\left({1-\overline{a}\cdot z}\right) ^{2n/p}}f\circ 
\psi (z)$,\ \par
is an isometry, i.e. $\left\Vert{Tf}\right\Vert _{H^{p}({\Bbb B})}=\left\Vert{f}\right\Vert 
_{H^{p}({\Bbb B})}$.\ \par
\end{Theorem}
\ \ \ \ \ \ Using it we get, exactly the same way as in the paper by Jevtic-Massaneda- 
Thomas~\cite{JevMasTho96},\ \par

\begin{Theorem} Let $\displaystyle S=\{a_{k}\}_{k\in {\Bbb N}}\subset {\Bbb B}$ 
be an interpolating sequence in $H^{p}({\Bbb B})$ and $\Phi $ be an automorphism 
of ${\Bbb B}$ then $\Phi (S):=\{\Phi (a_{k})\}_{k\in {\Bbb N}}\subset {\Bbb B}$ 
is still an interpolating sequence in $H^{p}({\Bbb B})$ with the same 
constant.\ \par
\end{Theorem}
\ \ \ \ \ \ In fact in~\cite{AmarThesis77} I explicitly stated a unitary 
representation of the group ${\mathrm{Aut}}({\Bbb B})$ in $H^{2}({\Bbb B})$ 
which gives as a corollary the fact that $H^{2}({\Bbb B})$ interpolating 
sequences are invariant by automorphisms.\ \par
\section{Interpolating sequences and discs.}
\ \ \ \ \ \ In this section we shall give the counter example in the unit 
ball ${\Bbb B}$ of ${\Bbb C}^{2}$  and for $1\leq p<\infty $.\ \par
\ \ \ \ \ \ If $S=\{a_{k},\ k\in {\Bbb N}\}$ is in $IH^{p}({\Bbb B})$ 
and if $(\sigma ,\phi )$ is a disc through $S$, $\displaystyle \sigma =\{\alpha _{k},\ k\in {\Bbb N}\}$, 
using automorphisms of the ball and of the disc, we may suppose that $a_{0}=0$ 
and $\alpha _{0}=0$ without changing the interpolating constants.\ \par
\ \ \ \ \ \ We have:\ \par
\ \ \ \ \ \ $\displaystyle \forall \lambda \in l_{H}^{p}(S),\ \exists f\in H^{p}({\Bbb 
B})\ s.t.\ \forall k\in {\Bbb N},\ f(a_{k})=\lambda _{k}$,\ \par
hence, for any $\displaystyle \lambda \in l_{H}^{p}(S)$, there is a $g\in A^{p}({\Bbb D})\ /\ g(\alpha _{k})=\lambda _{k}$, 
just take $g(z):=f\circ \phi (z)$.\ \par
So we have:\ \par

\begin{Proposition} If $S\in IH^{p}({\Bbb B})$ then any disc $(\sigma ,\phi )$ 
passing through $S$ is such that $\sigma $ is interpolating in $A^{p}({\Bbb D})$ 
for $l_{H}^{p}(S)$.\ \par
\end{Proposition}
\ \ \ \ \ \ But $\sigma $ in $IA^{p}({\Bbb D})$ means you can extend any 
$\lambda \in \ell _{A}^{p}(\sigma )$, and we have:\ \par

\begin{Lemma} If $(\sigma ,\phi )$ is a disc through $S$ with $a_{0}=0,\ \alpha _{0}=0$, 
we have that\ \par
\ $\forall k\in {\Bbb N},\ \left\vert{a_{k}}\right\vert \leq \left\vert{\alpha 
_{k}}\right\vert $ hence $1-\left\vert{a_{k}}\right\vert ^{2}\geq 1-\left\vert{\alpha _{k}}\right\vert 
^{2}$.\ \par
\end{Lemma}
\ \ \ \ \ \ Proof:\ \par
$\phi $ sends ${\Bbb D}$ into ${\Bbb B}$ hence it decreases the Gleason 
distance:\ \par
\ \ \ \ \ \ \ \ \ \ \ \ $d_{G}(a_{0},a_{k})\leq d_{G}(\alpha _{0},\alpha _{k})$;\ 
\par
because $d_{G}(a,b):={\mathrm{sup}}\{\left\vert{f(b)}\right\vert ;\ f\in H^{\infty 
},\ \left\Vert{f}\right\Vert _{\infty }\leq 1,\ f(a)=0\}$ and if $f\in H^{\infty }({\Bbb B})$ 
then\ \par
\ \ \ \ \ \ \ \ \ \ \ $g:=f\circ \phi \in H^{\infty }({\Bbb D})$ so\ \par
\ \ \ \ \ \ \ \ \ \ \ \ \ $d_{G}(\alpha ,\beta )={\mathrm{sup}}\{\left\vert{g(\beta )}\right\vert 
;\ g\in H^{\infty }({\Bbb D}),\ \left\Vert{g}\right\Vert _{\infty }\leq 
1,\ g(a)=0\}\geq d_{G}(a,b)$.\ \par
Because $a_{0}=0,\ \alpha _{0}=0$ we have $\displaystyle d_{G}(a_{0},a_{k})=\left\vert{a_{k}}\right\vert \leq \left\vert{\alpha 
_{k}}\right\vert =d_{G}(\alpha _{0},\alpha _{k})$ and the lemma.$\hfill 
\square$\ \par

\begin{Corollary} If $(\sigma ,\phi )\in {\cal{D}}_{S}$, $\lambda \in \ell _{H}^{p}(S){\Longrightarrow}\lambda \in \ell _{A}^{p}(\sigma 
)$ and $\left\Vert{\lambda }\right\Vert _{A,p}\leq \left\Vert{\lambda }\right\Vert 
_{H,p}$.\ \par
\end{Corollary}
\ \ \ \ \ \ Proof:\ \par
\ \ \ \ \ \ \ \ \ \ \ \ $\displaystyle \left\Vert{\lambda }\right\Vert _{H,p}^{p}:=\sum_{k\in 
{\Bbb N}}^{}{\left\vert{\lambda _{k}}\right\vert ^{p}(1-\left\vert{a_{k}}\right\vert 
^{2})^{2}}\geq \sum_{k\in {\Bbb N}}^{}{\left\vert{\lambda _{k}}\right\vert 
^{p}(1-\left\vert{\alpha _{k}}\right\vert ^{2})^{2}}=\left\Vert{\lambda 
}\right\Vert _{A,p}^{p}$.$\hfill \square$\ \par
\ \ \ \ \ \ Hence if $\sigma \subset {\Bbb D}$ is interpolating for $A^{p}({\Bbb D})$ 
it is interpolating in $A^{p}({\Bbb D})$ also for $\displaystyle \ell _{H}^{p}(S)$.\ 
\par
\ \ \ \ \ \ Now a natural question is:\ \par
\ \ \ \ \ \ \ \ \ \ \ \ if $S$ is a sequence in ${\Bbb B}$ containing 
$0$, such that any disc $(\sigma ,\phi )$ passing through $S$ is interpolating 
for $A^{p}({\Bbb D})$ then $S\in IH^{p}({\Bbb B})$  ?\ \par
\ \ \ \ \ \ We notice that this is stronger, by the previous corollary, 
that the necessary  condition which is:\ \par
\ \ $\sigma $  is interpolating in $A^{p}({\Bbb D})$  for $\ell _{H}^{p}(S)$.\ 
\par
\ \ \ \ \ \ Despite this stronger condition the answer is no:\ \par

\begin{Theorem} For any $p\in [1,\ \infty [$, there is a sequence $S$ 
in ${\Bbb B}$ such that all disc passing through $S$ is interpolating 
for $A^{p}({\Bbb D})$ with a uniformly bounded constant but $S$ is not 
in $IH^{p}({\Bbb B})$.\ \par
\end{Theorem}
\ \ \ \ \ \ Proof:\ \par
take $1\leq p<\infty $, we know that if a sequence $\sigma $ in ${\Bbb D}$ 
is enough separated, i.e.\ \par
\ \ \ \ \ \ \ $\forall (\alpha ,\beta )\in S^{2}\ \alpha {\not =}\beta ,\ d_{G}(\alpha 
,\beta )\geq \delta $,\ \par
for a certain $\delta =\delta _{p}>0$ then $\sigma $ is in $IA^{p}({\Bbb D})$ 
with an interpolating constant depending only on $\delta $ ~\cite{amBerg78}. 
In~\cite{Seip93}, K. Seip characterized the interpolating sequences for 
$A^{p}({\Bbb D})$ but here it suffices to use the result in~\cite{amBerg78} 
which is true also in several variables.\ \par
\ \ \ \ \ \ Take a sequence $S$ in ${\Bbb B}$ such that $S$ cannot be 
contained in the zero set of a $H^{p}({\Bbb B})$ function, but such that 
the points in $S$  are separated by $\delta $.\ \par
\ \ \ \ \ \ This is easy to do: as in ~\cite{amBerg78} just take for $S$ 
a full net of points in ${\Bbb B}$  separated by $\delta $. The following 
lemma tells us that $S$ cannot be the zero set of any $H^{p}({\Bbb B})$ 
 function.\ \par

\begin{Lemma} Suppose that the sequence $S$ is such that for any point 
$\zeta \in \partial {\Bbb B}$ the admissible convergence region $\Gamma (\zeta ,\alpha ):=\{z\in {\Bbb B}\ /\ \left\vert{1-\overline{\zeta 
}\cdot z}\right\vert <\alpha (1-\left\vert{z}\right\vert ^{2})$ of aperture 
$\alpha $ contains an infinite number of points of $S$, then $S$ cannot 
be contained in a non-trivial $H^{p}({\Bbb B})$ zero set.\ \par
\end{Lemma}
\ \ \ \ \ \ Proof of the lemma:\ \par
suppose that there is an $H^{p}({\Bbb B})$ function $f$ which is zero 
on $S$. We know that the  Fatou convergence theorem is valid, which means 
that $f$ admits admissible boundary values $f^{*}$ almost everywhere on 
$\partial {\Bbb B}$ and $f$ is the Poisson integral of $f^{*}$ (~\cite{RudinBall81} 
\S 5.4, p72).\ \par
\ \ \ \ \ \ Hence we have $f^{*}(\zeta )={\mathrm{li}}{\mathrm{m}}_{a\in S\cap \Gamma (\zeta ,\alpha 
),\ a{\longrightarrow}\zeta }f(a),\ a.e.\ \zeta \in \partial {\Bbb B}$; 
thus $f^{*}=0\ a.e.$ hence $f=0$. $\hfill \square$\ \par
\ \ \ \ \ \ Taking  a full net of points in ${\Bbb B}$ separated by $\delta $, 
one sees easily that, with $\alpha $  big enough with respect to $\delta $, 
the condition of the lemma is satisfied.\ \par
\ \ \ \ \ \ Then for any disc $(\sigma ,\phi )$  passing through $S$, 
the points of $\sigma $ are separated by at least $\delta $, hence $\sigma \in IA^{p}({\Bbb D})$ 
 with a uniform constant.\ \par
\ \ \ \ \ \ But $S$ cannot be in $IH^{p}({\Bbb B})$ because $S$ is not 
contained in a non-trivial $H^{p}({\Bbb B})$  zero set and this is a necessary 
condition:\ \par
\ \ \ \ \ \ \ if $S$ is $H^{p}({\Bbb B})$ interpolating, then there is 
a function $f\in H^{p}({\Bbb B})$ such that $f(a_{1})=1,\ \forall k>1,\ f(a_{k})=0$. 
Hence $f$ is zero on $S\backslash \{a_{1}\}$ and of course $f$ is not 
identically $0$. So the function $\displaystyle g(z):=\frac{\overline{a}_{1}\cdot (z-a_{1})}{1-\overline{a}_{1}\cdot 
z}f(z)$ is still in $H^{p}({\Bbb B})$, still not identically $0$, and 
is zero on $S$, which proves that $S$ must be contained in the zero set 
of a $H^{p}({\Bbb B})$  function.$\hfill \square$\ \par

\bibliographystyle{/usr/share/texmf/bibtex/bst/base/plain}
\bibliography{/home/eamar/bibFiles/General.bib,/home/eamar/bibFiles/Amar.bib,/home/eamar/bibFiles/Interpol.bib}

\begin{thebibliography}{10}

\bibitem{AmarThesis77}
E.~Amar.
\newblock {\em Suites D'interpolation Dans Le Spectre D'une Alg\`ebre
  D'op\'erateurs}.
\newblock PhD thesis, Universit\'e Paris XI, Orsay, 1977.

\bibitem{amBerg78}
E.~Amar.
\newblock Suites d'interpolation pour les classes de {B}ergman de la boule et
  du polydisque de $c^n$.
\newblock {\em Canadian J. Math.}, 30:711--737, 1978.

\bibitem{amCor91}
E.~Amar.
\newblock On the corona problem.
\newblock {\em The Journal of Geometric Analysis}, 1(4):291--305, 1991.

\bibitem{AmMenOp02}
E.~Amar and C.~Menini.
\newblock On an operator theory approach to the corona problem.
\newblock {\em Bull. London Math. Soc.}, 34:369--373, 2002.

\bibitem{AmTho94}
E.~Amar and P.~Thomas.
\newblock A notion of extremal discs related to interpolation in the ball.
\newblock {\em Math. Ann.}, 300:419--433, 1994.

\bibitem{Duren70}
P.~Duren.
\newblock {\em Theory of ${H}^p$ spaces}, volume~38 of {\em Pure and Applied
  Mathematics}.
\newblock Academic Press, New-York, 1970.

\bibitem{GlobevnikDiscs88}
J.~Globevnik.
\newblock Discs in the ball containing given discrete sets.
\newblock {\em Math. Ann.}, 281:87--96, 1988.

\bibitem{JevMasTho96}
Miroljub Jevtic, Xavier Massaneda, and Pascal~J. Thomas.
\newblock Interpolating sequences for weighted {B}ergman spaces of the ball.
\newblock {\em Mich. Math. J. 43, No.3, 495-517 (1996)}, 43(3):495--517, 1996.

\bibitem{RudinBall81}
W.~Rudin.
\newblock Function theory in the unit ball of $ {C}^n$.
\newblock {\em Grundenlehren}, 1981.

\bibitem{Seip93}
K.~Seip.
\newblock Beurling type density theorems in the unit disk.
\newblock {\em Inventiones mathmatic\ae}, 113:21--39, 1993.

\end{thebibliography}

\end{document}